\def\MS{\noalign{\medskip}}

\def\midline#1{\setbox1=\hbox{\vrule height 0.5pt depth 0pt width #1}
              \raise 1.60pt \box1}
\def\down#1{\vphantom{\hbox{\vrule height #1pt depth 0pt width 0pt}}}
\def\up#1{\vphantom{\lower #1pt\hbox{\vrule height #1pt depth #1pt width 0pt}}}
 \newcount\refnum
\refnum =0
\def\ref{\advance\refnum by 1\smallbreak\item{\the\refnum.}} \def\refsq{\advance\refnum by 1\smallbreak\item{[\the\refnum]}} 
%
\def\buildrel#1\over#2{\mathrel{
\mathop{\kern 0pt#2}\limits^{#1}}}
\def\bbuildrel#1_#2^#3{\mathrel{
\mathop{\kern 0pt#1}\limits_{#2}^{#3}}}

\def\newpage{\vfill\break}

\def\newline{}

\def\ts{{\textstyle{1\over 2}}}
\font\title =cmbx10 at 12pt

\font\author =cmr10 at 12pt
\font\abstract =cmr9
\font\sevenbf =cmbx7

\magnification =1200
\parindent=0pt
\headline{\hfill \the\pageno}
\footline{\hfil}
\null\vskip 1truein
%
%
\centerline{{\title ANALYTIC CONTINUATION OF THE GENERALIZED}} \centerline{{\title HYPERGEOMETRIC SERIES NEAR UNIT ARGUMENT }} \centerline{{\title WITH EMPHASIS ON THE ZERO-BALANCED SERIES}} \bigskip
\bigskip
\centerline{{\bf WOLFGANG B\"UHRING}}
\smallskip
\centerline{Physikalisches Institut, Universit\"at Heidelberg} \centerline{69120 Heidelberg, Germany} \centerline{{\sevenbf E-Mail: BUEHRING@PHYSI.UNI-HEIDELBERG.DE}} \medskip
\centerline{and}
\medskip
\centerline{{\bf H.M. SRIVASTAVA}}
\smallskip
\centerline{Department of Mathematics and Statistics, 
University of Victoria}
\centerline{Victoria, British Columbia~~V8W 3P4, Canada} \centerline{{\sevenbf E-Mail: HMSRI@UVVM.UVIC.CA}}
\footnote{}{2000 {\it Mathematics Subject Classification}. 
Primary 33C20, 34E05; Secondary 41A58.}
\bigskip
\bigskip
\centerline{{\bf ABSTRACT}}
\medskip
{\parindent=0.5truein\narrower
Various methods to obtain the analytic continuation near $z=1$ of the hypergeometric
series $_{p+1}F_p(z)$ are reviewed together with some of the results. One approach is to establish a recurrence relation with respect to $p$ and then, after its repeated use, to resort to the well-known properties of the Gaussian hypergeometric series. Another approach is based on the properties of the underlying generalized hypergeometric differential equation: For the coefficients in the power series expansion around $z=1$ a general formula, valid for any $p$, is found in terms of a limit of a terminating Saalsch\"utzian hypergeometric series of unit argument. Final results may then be obtained for each particular $p$ after application of an appropriate transformation formula of the Saalsch\"utzian hypergeometric series. The behaviour at unit argument of zero-balanced hypergeometric series, which have received particular attention in recent years, is discussed in more detail. The related problem involving the behaviour of partial sums of such series is addressed briefly.\par}
\bigskip
\bigskip
\baselineskip=14pt
\parindent =0.5truein
\newpage
\noindent{\bf 1.~~Introduction}
\medskip
Generalized hypergeometric series like the Gaussian hypergeometric series but with more parameters have received renewed attention in recent years, especially their value (if finite) at or behaviour near unit argument. This interest was partly motivated by some result of this type contained without proof in the notebooks of Ramanujan ({\it cf}. [15], [3], [10], and [4]), partly by the need for such information in the context of generalized hypergeometric functions of several variables ({\it cf}. [16] and [20]). In this article we are concerned with hypergeometric series or functions (see, {\it e.g.}, [2] and [18]): $$\eqalign{&_{p+1}F_p\left(\matrix{\up{4}\hfill a_1,a_2,\cdots ,a_{ p+1};\cr
\down{4}\hfill b_1,b_2,\cdots ,b_p;\cr}
\,\,\,z\right)\cr
\MS&\qquad\qquad\qquad =\sum_{n=0}^{\infty}\,{{(a_1)_n(a_2)_n\cdots (a_{p+1})_n}\over {(b_1)_n(b_2)_n\cdots (b_p)_n}}\,\,{{z^n}\over { n!}}\qquad (|z|<1),\cr}
\eqno (1.1)$$
written in terms of Pochhammer symbols $(\lambda )_n$, where $$(\lambda )_n=\lambda (\lambda +1)\cdots (\lambda +n-1)= {{\Gamma (\lambda +n)}\over {\Gamma (\lambda )}}.$$ An important quantity for such functions is $$s=s_p=\sum_{j=1}^p\,b_j-\sum_{j=1}^{p+1}\,a_j\,,\eqno (1.2)$$ the sum of the denominator parameters minus the sum of the numerator parameters, which is the non-trivial characteristic exponent of the underlying (generalized) hypergeometric differential equation at unit argument and so determines the behaviour of the hypergeometric function near this point. We may also recall that a hypergeometric series (1.1) with $s$ equal to an integer is called $s$-balanced and that a one-balanced series is called a Saalsch\"utzian series.

\headline{\hfil -\the\pageno -\hfil}

It is the zero-balanced case which has received particlar attention, beginning with the following formula in the notebook of Ramanujan ([15], [3], [10], and [4]) for parameters such that $a_1+a_2+a_3=b_1+b_2$:
$$\eqalign{&{{\Gamma (a_1)\Gamma (a_2)\Gamma (a_3)}\over {\Gamma ( b_1)\Gamma (b_2)}}\,\,\,_3F_2\left(\matrix{\up{4}\hfill a_1,a_2,a_ 3;\cr
\down{4}\hfill b_1,b_2;\cr}
\,\,\,z\right)\cr
\MS&\qquad\qquad\qquad =-\log (1-z)+L+O\left((1-z)\log (1-z)\right ),\qquad z\to 1,\cr}
\eqno (1.3)$$
with, if $\Re (a_3)>0$,
$$L=2\psi (1)-\psi (a_1)-\psi (a_2)+\sum_{k=1}^{\infty}\,{{(b_2-a_ 3)_k(b_1-a_3)_k}\over {k(a_1)_k(a_2)_k}},\eqno (1.4)$$ where
$$\psi (z)={{\Gamma^{\prime}(z)}\over {\Gamma (z)}}$$ is the logarithmic derivative of the gamma function. Early proofs of (1.3) and (1.4) were supplied by Saigo [16] and Berndt [3], and by Evans and Stanton [9]. In a more general setting, the subject was reconsidered by Saigo and H.M. Srivastava [17] and by B\"uhring ([5] and [6]). Finally, A.K. Srivastava [19] addressed the related but different problem, earlier treated by Ramanujan [15] as well as by Evans and Stanton [9], of the asymptotic behaviour of partial sums of zero-balanced hypergeometric series.

In the present article we shall portray the different methods available in the literature, and review some of the results obtained by these methods.
\bigskip
\bigskip\newline
\noindent{\bf 2.~~Recurrence with respect to $p$} \bigskip
One approach is to establish a recurrence relation with respect to $p$ between generalized hypergeometric series (1.1) and then to apply it repeatedly until one arrives at the Gaussian hypergeometric series, for which all the required properties are known. In particular, it may be helpful to remember that, if
$$s_1=b_1-a_1-a_2\eqno (2.1)$$
is not equal to an integer, the analytic continuation near $z=1$ is given by the continuation formula:
$$\eqalign{&{{\Gamma (a_1)\Gamma (a_2)}\over {\Gamma (b_1)}}\,\,_2 F_1\left(\matrix{\up{4}\hfill a_1,a_2;\cr \down{4}\hfill b_1;\cr}
\,\,z\right)=\sum_{n=0}^{\infty}(-1)^n\,{{\Gamma (a_1+n)\Gamma (a_ 2+n)\Gamma (s_1-n)}\over {\Gamma (a_1+s_1)\Gamma (a_2+s_1)\,}}\,\,{{ (1-z)^n}\over {n!}}\cr
\MS&\qquad +(1-z)^{s_1}\,\sum_{n=0}^{\infty}(-1)^n\,{{\Gamma (a_1+ s_1+n)\Gamma (a_2+s_1+n)\Gamma (-s_1-n)}\over {\Gamma (a_1+s_1)\Gamma (a_2+s_1)}}\,\,{{(1-z)^n}\over {n!}}.\cr} \eqno (2.2)$$
Also, if the real part of $s_1$ is positive, the Gaussian hypergeometric function is finite at $z=1$ and its value is given by the Gaussian summation formula:
$${1\over {\Gamma (b_1)}}\,\,_2F_1\left(\matrix{\hfill a_1,a_2;\cr \hfill b_1^{};\cr}
\,\,1\right)={{\Gamma (b_1-a_1-a_2)}\over {\Gamma (b_1-a_1)\Gamma (b_1-a_2)}}={{\Gamma (s_1)}\over {\Gamma (a_1+s_1)\Gamma (a_2+s_1)}} ,\eqno (2.3)$$
and this formula may be viewed as a limiting case of the formula of Saalsch\"utz:
$${1\over {\Gamma (b_1)}}\,\,_3F_2\left(\matrix{\up{4}\hfill a_1,a_ 2,-m;\cr
\down{4}\hfill b_1,1-s_1-m;\cr}
\,\,1\right)={{(a_1+s_1)_m(a_2+s_1)_m}\over {(s_1)_m\,\Gamma (b_1+ m)}},\eqno (2.4)$$
which is valid for $m=0,1,2,\cdots ,$ if $s_1$ is not equal to a negative integer or zero. Finally, we remember that, when $z\to 1$, the zero-balanced series has the behaviour:
$${{\Gamma (a_1)\Gamma (a_2)}\over {\Gamma (a_1+a_2)}}\,\,_2F_1\left (\matrix{\up{4}\hfill a_1,a_2;\cr
\down{4}\hfill a_1+a_2;\cr}
\,\,z\right)=2\psi (1)-\psi (a_1)-\psi (a_2)-\log (1-z)+o(1).\eqno (2.5)$$

Assuming that $\Re (a_{p+1})>0$, we may use the Gaussian summation formula to get
$$\eqalign{&{{\Gamma (a_1)\Gamma (a_2)\cdots\Gamma (a_{p+1})}\over { \Gamma (b_1)\cdots\Gamma (b_p)}}\,\,\,_{p+1}F_p\left(\matrix{\up{4}\hfill a_1,a_2,\cdots ,a_{p+1};\cr
\down{4}\hfill b_1,\cdots ,b_p;\cr}
\,\,z\right)\cr
\MS&\qquad =\sum_{n=0}^{\infty}\,{{\Gamma (a_1+n)\cdots\Gamma (a_p +n)\Gamma (a_{p+1}+n)}\over {\Gamma (b_1+n)\cdots\Gamma (b_{p-2}+n )\Gamma (b_{p-1}+n)\Gamma (b_p+n)}}\,\,{{z^n}\over {n!}}\cr \MS&\qquad =\sum_{n=0}^{\infty}\,{{\Gamma (a_1+n)\cdots\Gamma (a_p +n)}\over {\Gamma (b_1+n)\cdots\Gamma (b_{p-2}+n)}}\cr \MS&\qquad\qquad\cdot{1\over {\Gamma (b_p+b_{p-1}-a_{p+1}+n)}}\,\,_ 2F_1\left(\matrix{\up{4}\hfill b_p-a_{p+1},b_{p-1}-a_{p+1};\cr \down{4}\hfill b_p+b_{p-1}-a_{p+1}+n;\cr} \,\,\,1\right)\,{{z^n}\over {n!}}.\cr}
\eqno (2.6)$$
Writing the $_2F_1(1)$ as a series and interchanging the order of the summations, we may obtain the required recurrence relation ([17] and [6]):
$$\eqalign{&_{p+1}F_p\left(\matrix{\up{4}\hfill a_1,a_2,\cdots ,a_{p+1};\cr
\down{4}\hfill b_1,\cdots ,b_p;\cr}
\,\,z\right)={{\Gamma (b_p)\Gamma (b_{p-1})}\over {\Gamma (a_{p+1} )\Gamma (b_p+b_{p-1}-a_{p+1})}}\cr &  \quad\cdot\sum_{m=0}^{\infty}\,{{(b_p-a_{p+1})_m(b_{p-1}-a_{p+ 1})_m}\over {(b_p+b_{p-1}-a_{p+1})_m\,m!}}\,\,_pF_{p-1}\left(\matrix{\up{ 4}\hfill a_1,a_2,\cdots ,a_p;\cr
\down{4}\hfill b_1,\cdots ,b_{p-2},b_p+b_{p-1}-a_{p+1}+m;\cr} \,\,z\right),\cr
}\eqno(2.7)$$
valid if $\Re (a_{p+1})>0.$ This recurrence relation may be used as it stands [6] or first be rewritten [17] as follows with the term for $m=0$ separated:
$$\eqalign{&_{p+1}F_p\left(\matrix{\up{4}\hfill a_1,a_2,\cdots ,a_{ p+1};\cr
\down{4}\hfill b_1,\cdots ,b_p;\cr}
\,z\right)\cr
\MS&\quad ={{\Gamma (b_p)\Gamma (b_{p-1})}\over {\Gamma (a_{p+1})\Gamma (b_p+b_{p-1}-a_{p+1})}}\,\,_pF_{p-1}\left(\matrix{\up{4}\hfill a_1 ,a_2,\cdots ,a_p;\cr
\down{4}\hfill b_1,\cdots ,b_{p-2},b_p+b_{p-1}-a_{p+1};\cr} \,z\right)\cr
\MS&\qquad\qquad +{{(b_p-a_{p+1})(b_{p-1}-a_{p+1})\Gamma (b_p)\Gamma (b_{p-1})}\over {\Gamma (a_{p+1})\Gamma (b_p+b_{p-1}-a_{p+1}+1)}}\cr \MS&\qquad\cdot\sum_{m=0}^{\infty}\,{{(b_p-a_{p+1}+1)_m(b_{p-1}-a_{ p+1}+1)_m}\over {(b_p+b_{p-1}-a_{p+1}+1)_m(m+1)!}}\cr \MS&\qquad\cdot\,_pF_{p-1}\left(\matrix{\up{4}\hfill a_1,a_2,\cdots ,a_p;\cr
\down{4}\hfill b_1,\cdots ,b_{p-2},b_p+b_{p-1}-a_{p+1}+m+1;\cr} \,z\right).\cr}
\eqno (2.8)$$
In both cases, the sum over $m$ may be viewed as a special Kamp\'e de F\'eriet series, ({\it cf}., {\it e.g.}, [20]), but we do not want to introduce and use the notation for such series here. While (2.7) can be applied in any case, (2.8) is specially designed for the case when the series on the left is zero-balanced. Then only the first term on the right is singular at $z=1$, but the hypergeometric series in the second term is $(m+1)$-balanced and therefore finite at $z=1$ and can be summed by the Gaussian summation formula (2.3) when $p=2$ (or by its generalization when $p>2$). The results so obtained [17], that is, the leading terms of the behaviour when $z\to 1$ of the zero-balanced hypergeometric series, shall be displayed below in Section~4. 

Alternatively, using (2.7) repeatedly and rearranging the summations, we finally arrive at [6]
$$\eqalign{&{{\Gamma (a_1)\Gamma (a_2)\cdots\Gamma (a_{p+1})}\over { \Gamma (b_1)\cdots\Gamma (b_p)}}\,\,_{p+1}F_p\left(\matrix{\up{4}\hfill a_1,a_2,\cdots ,a_{p+1};\cr
\down{4}\hfill b_1,\cdots ,b_p;\cr}
\,z\right)\cr
\MS&\qquad =\sum_{k=0}^{\infty}\,A_k^{(p)}\,{{\Gamma (a_1)\Gamma ( a_2)}\over {\Gamma (s+a_1+a_2+k)}}\,\,_2F_1\left(\matrix{\up{4}\hfill a_1,a_2;\cr
\down{4}\hfill s+a_1+a_2+k;\cr}
\,z\right),\cr}
\eqno (2.9)$$
where
$$A_k^{(2)}={{(b_2-a_3)_k(b_1-a_3)_k}\over {k!}},\eqno (2.10)$$ $$\eqalign{
A_k^{(3)}&={{(b_3+b_2-a_4-a_3)_k(b_1-a_3)_k}\over {k!}}\cr
\MS&\qquad\cdot\, _3F_2\left(\matrix{
\up{4}\hfill b_3-a_4,b_2-a_4,-k;\cr
\down{4}\hfill b_3+b_2-a_4-a_3,1+a_3-b_1-k;\cr} \,1\right),}\eqno (2.11)$$
$$\eqalign{&A_k^{(4)}={{(b_4+b_3+b_2-a_5-a_4-a_3)_k(b_1-a_3)_k}\over { k!}}\cr
\MS&\qquad\qquad\cdot\sum_{\ell =0}^k\,{{(b_4+b_3-a_5-a_4)_{\ell}( b_2-a_4)_{\ell}(-k)_{\ell}}\over {(b_4+b_3+b_2-a_5-a_4-a_3)_{\ell} (1+a_3-b_1-k)_{\ell}\,\ell !}}\cr
\MS&\qquad\qquad\cdot\,_3F_2\left(\matrix{\up{4}\hfill b_4-a_5,b_3 -a_5,-\ell ;\cr
\down{4}\hfill b_4+b_3-a_5-a_4,1+a_4-b_2-\ell ;\cr} \,1\right),\cr}
\eqno (2.12)$$
and
$$\eqalign{
&A_k^{(p)}={{(b_p+b_{p-1}+\cdots +b_2-a_{p+1}-a_p-\cdots -a_3)_k(b_1-a_3)_k}\over {k!}}\cr
\MS&\cdot\sum_{k_2=0}^k\,{{(-k)_{k_2}}\over {(b_p+b_{ p-1}+\cdots +b_2-a_{p+1}-a_p-\cdots -a_3)_{k_2}(1+a_3-b_1-k)_{k_2}}}\cr \MS&\cdot{{(b_p+b_{p-1}+\cdots +b_3-a_{p+1}-a_p-\cdots -a_4)_{k_2}(b_2-a_4)_{k_2}}\over {k_2!}}\cr \MS&\cdot\sum_{k_3=0}^{k_2}\,{{(-k_2)_{k_3}}\over {(b_ p+b_{p-1}+\cdots +b_3-a_{p+1}-a_p-\cdots -a_4)_{k_3}(1+a_4-b_2-k_2 )_{k_3}}}\cr
\MS&\cdot\,\cdots\,\cdot {{(b_p+b_{p-1}-a_{p+1}-a_p)_{k_{p-2}}( b_{p-2}-a_p)_{k_{p-2}}}\over {k_{p-2}!}}\cr \MS&\cdot\sum_{k_{p-1}=0}^{k_{p-2}}\,{{(-k_{p-2})_{k_{ p-1}}}\over {(b_p+b_{p-1}-a_{p+1}-a_p)_{k_{p-1}}(1+a_p-b_{p-2}-k_{ p-2})_{k_{p-1}}}}\cr
\MS&\qquad\qquad\qquad\cdot {{(b_p-a_{p+1})_{k_{p-1}} (b_{p-1}-a_{p+1})_{k_{p-1}}}\over {k_{p-1}!}},\cr} \eqno (2.13)$$
provided that $\Re (a_j)>0$ for $j=3,4,\cdots ,p+1$. 

Substituting for the $_2F_1(z)$ in (2.9) its analytic continuation near $z=1$, we get, if $s$ is not equal to an integer, for $|1-z|<1$, $|\arg (1-z)|<\pi$, and $p=2,3,\cdots ,$ $$\eqalign{&{{\Gamma (a_1)\Gamma (a_2)\cdots\Gamma (a_{p+1})}\over { \Gamma (b_1)\cdots\Gamma (b_p)}}\,\,_{p+1}F_p\left(\matrix{\up{4}\hfill a_1,a_2,\cdots ,a_{p+1};\cr
\down{4}\hfill b_1,\cdots ,b_p;\cr}
\,z\right)\cr
\MS&\qquad =\sum_{n=0}^{\infty}\,g_n(0)(1-z)^n+(1-z)^s\,\sum_{n=0}^{ \infty}\,g_n(s)(1-z)^n\cr}
\eqno (2.14)$$
with
$$\eqalign{&g_n(r)=(-1)^n\,\,{{\Gamma (a_1+r+n)\Gamma (a_2+r+n)\Gamma (s-2r-n)}\over {\Gamma (a_1+s)\Gamma (a_2+s)\,n!}}\cr \MS&\qquad\qquad\qquad\cdot\sum_{k=0}^{\infty}\,\,{{(s-r-n)_k}\over { (a_1+s)_k(a_2+s)_k}}\,\,A_k^{(p)}.\cr}
\eqno (2.15)$$
Here the series terminates when $r=s$, but in the case when $r=0$, the condition:
$$\Re (a_3+n)>0\land\cdots\land\Re (a_{p+1}+n)>0$$ is required for convergence.

A related result, which is obtained from (2.9) with $z=1$, is a generalization of the Gaussian summation formula, for $\Re (s)>0$, and we have
$$\eqalign{&{{\Gamma (a_1)\Gamma (a_2)\cdots\Gamma (a_{p+1})}\over { \Gamma (b_1)\cdots\Gamma (b_p)}}\,\,_{p+1}F_p\left(\matrix{\up{4}\hfill a_1,a_2,\cdots ,a_{p+1};\cr
\down{4}\hfill b_1,\cdots ,b_p;\cr}
\,1\right)\cr
\MS&\qquad ={{\Gamma (a_1)\Gamma (a_2)\Gamma (s)}\over {\Gamma (a_ 1+s)\Gamma (a_2+s)}}\,\,\sum_{k=0}^{\infty}\,{{(s)_k}\over {(a_1+s )_k(a_2+s)_k}}\,\,A_k^{(p)},\cr}
\eqno (2.16)$$
where the series converges if
$$\Re (a_3)>0\land\cdots\land\Re (a_{p+1})>0.$$ For $p=2$, (2.16) reduces to a known formula ([2] and [11]). \bigskip
\bigskip\newline
\noindent{\bf 3.~~Limit formulas}
\bigskip
Another approach is to get the coefficients in the continuation formula as certain \penalty-10000 limits of terminating Saalsch\"utzian hyper\-geometric
series. The generalized hypergeometric function or series (1.1) is a particular solution of a $(p+1)$-th order linear differential equation with three regular singular points at $z=0,1,\infty$. The characteristic exponents of this differential equation at $z=1$ are $0,1,2,\cdots ,p-1$ and $s$
defined by (1.2). If $s$ is not equal to an integer, then all the exponents (despite the occurrence of integral differences between them) give solutions of the differential equation which are linearly independent of each other, and therefore the analytic continuation near $z=1$ of the generalized hyper\-geometric function is then given by (2.14), already obtained above by a different reasoning. While it is also known from the detailed investigation by N\o rlund ({\it cf}. [13] and [12]) that the connection coefficient $g_0(s)$ is equal to $\Gamma (-s)$ for any $ p$, the
remaining connection coefficients $g_0(0),g_1(0),\cdots ,g_{p-1}(0 )$ could only more
recently be determined, first in [5] for $p=2$ and later in [6], using the method of Section~2, for any $p$. 

The coefficient of the leading singular contribution to (2.14) (at $ z=1$)
can most conveniently be obtained by Darboux's method [14] as described in [5] for $p=2$. Applied to any $p$, this method yields the same result:
$$g_0(s)=\Gamma (-s),\eqno (3.1)$$
if presented in terms of $s$, which according to (1.2) depends on $ p$.

In order to obtain $g_0(0)$ by the same method, we multiply (2.14) by $(1-z)^{-s}$. The left-hand side then becomes $$\eqalign{&\sum_{n=0}^{\infty}\,{{(s)_n}\over {n!}}\,z^n\,\,\sum_{ k=0}^{\infty}\,\,{{\prod\limits_{j=1}^{p+1}\,\Gamma (a_j+k)}\over { \prod\limits_{j=1}^p\,\Gamma (b_j+k)}}\,\,{{z^k}\over {k!}}\cr \MS&\qquad =\sum_{m=0}^{\infty}\,\,\sum_{k=0}^m\,\,{{\prod\limits_{ j=1}^{p+1}\,\Gamma (a_j+k)}\over {\prod\limits_{j=1}^p\,\Gamma (b_ j+k)}}\,{{(s)_{m-k}}\over {k!(m-k)!}}\,z^m.\cr} \eqno (3.2)$$
If this is written in the form $\Sigma u_mz^m$, the coefficients $ u_m$, after
application of the identity:
$${{(s)_{m-k}}\over {(m-k)!}}={{(s)_m(-m)_k}\over {m!(1-s-m)_k}},\eqno (3.3)$$
are
$$u_m={{\prod\limits_{j=1}^{p+1}\,\Gamma (a_j)}\over {\prod\limits_{ j=1}^p\,\Gamma (b_j)}}\,\,{{(s)_m}\over {m!}}\,\,\,_{p+2}F_{p+1}\left (\matrix{\up{4}\hfill a_1,a_2,\cdots ,a_{p+1},-m;\cr \down{4}\hfill b_1,b_2,\cdots ,b_p,1-s-m;\cr} \,1\right).\eqno (3.4)$$
On the right-hand side the leading singular term is $g_0(0)(1-z)^{ -s}$,
which, when expanded in powers of $z$ and written in the form $\Sigma w_mz^m$, has coefficients
$$w_m=g_0(0)\,\,{{(s)_m}\over {m!}}.\eqno (3.5)$$ Comparison of the $u_m$ and $w_m$, which have to agree asymptotically as $m\to\infty$, then yields
$$g_0(0)={{\prod\limits_{j=1}^{p+1}\,\Gamma (a_j)}\over {\prod\limits_{ j=1}^p\,\Gamma (b_j)}}\,\,\lim_{m\to\infty}\,\,_{p+2}F_{p+1}\left(\matrix{\up{ 4}\hfill a_1,a_2,\cdots ,a_{p+1},-m;\cr
\down{4}\hfill b_1,b_2,\cdots ,b_p,1-s-m;\cr} \,1\right).\eqno (3.6)$$
In order to obtain the $g_n(0)$ for $n>0$ too, we consider the $n$-th derivative with respect to $z$ of (2.14) multiplied by $(1-z)^{-s+ n}$, so that
the term with $g_n(0)$ becomes the leading singular term. Proceeding as above, we get, as a generalization of (3.6), $$\eqalign{&g_n(0)={{(-1)^n}\over {n!}}\,\,{{\prod\limits_{j=1}^{p +1}\,\Gamma (a_j+n)}\over {\prod\limits_{j=1}^p\,\Gamma (b_j+n)}}\cr \MS&\qquad\qquad\quad\cdot\lim_{m\to\infty}\,\,_{p+2}F_{p+1}\left(\matrix{\up{ 4}\hfill a_1+n,a_2+n,\cdots ,a_{p+1}+n,-m;\cr \down{4}\hfill b_1+n,b_2+n,\cdots ,b_p+n,1-s+n-m;\cr} \,1\right).\cr}
\eqno (3.7)$$
The $g_n(s)$ for $n>0$ could be derived from $g_0(s)$ and the recurrence relation which they satisfy [13], but it is desirable to get an explicit representation the following way. The \penalty-10000 coefficients $g_n(r)$ satisfy a recurrence relation which is obtained from the differential equation by inserting, according to the Frobenius method, a power-series solution:
$$\sum_{n=0}^{\infty}\,g_n(r)(1-z)^{r+n}.$$ The recurrence relation, therefore, depends on the exponent $r$ {\it via} the combination $r+n$ only. Starting from the $g_n(0)$ of (3.7), rewriting the factor
$${{(-1)^n}\over {n!}}={{\sin\left(\pi (s-n)\right)}\over {\Gamma (1+n)\sin (\pi s)}}={{\Gamma (1-s)\Gamma (s)}\over {\Gamma (1+n)\Gamma (1-s+n)\Gamma (s-n)}},\eqno (3.8)$$
and replacing $n$ by $r+n$, we obtain a solution of the recurrence relation corresponding to the exponent $r$. For $r=s$, this solution can differ from the desired $g_n(s)$ only by a constant factor which does not depend on $n$. We may, therefore, write $$\eqalign{&g_n(r)={{V(r)\Gamma (1-s)\Gamma (s)}\over {\Gamma (1+r +n)\Gamma (1-s+r+n)}}\,\,{{\prod\limits_{j=1}^{p+1}\,\Gamma (a_j+r +n)}\over {\prod\limits_{j=1}^p\,\Gamma (b_j+r+n)}}\,\,\lim_{m\to\infty}\,\,{ 1\over {\Gamma (s-r-n)}}\cr
\MS&\qquad\cdot\,_{p+2}F_{p+1}\left(\matrix{\up{4}\hfill a_1 +r+n,a_2+r+n,\cdots ,a_{p+1}+r+n,-m;\cr
\down{4}\hfill b_1+r+n,b_2+r+n,\cdots ,b_p+r+n,1-s+r+n-m;\cr} \,1\right),\cr}
\eqno (3.9)$$
where $V(0)=1$ and $V(s)$ is such that (3.1) holds true. In order to determine $V(s)$, we first consider the case when $p=1$. Then the hypergeometric series in (3.9) can be summed by means of Saalsch\"utz's formula (2.4), which yields
$${{(a_1+s)_m(a_2+s)_m}\over {(s-r-n)_m(b_1+r+n)_m}}\quad 
\bbuildrel\hbox to .4truein{\rightarrowfill}_{m\to\infty}^{}\quad {{\Gamma (s-r-n)\Gamma (b_1+r+n)}\over {\Gamma (a_1+s)\Gamma (a_2+s)}} .\eqno (3.10)$$
The factor $\Gamma (s-r-n)$, which is singular for $r=s$, now cancels in (3.9) and we have, for $p=1$,
$$g_n(r)={{V(r)\Gamma (1-s)\Gamma (s)}\over {\Gamma (1+r+n)\Gamma (1-s+r+n)}}\,\,{{\Gamma (a_1+r+n)\Gamma (a_2+r+n)}\over {\Gamma (a_ 1+s)\Gamma (a_2+s)}}.\eqno (3.11)$$
For $r=s$ and $n=0$, comparison with (3.1) shows that $V(s)=-1$. This value is inde\-pendent of the parameters of the hypergeometric function and so it holds true for any $p$. The first fraction on the right of (3.9) may now be rewritten as follows:
$$\eqalign{{{V(r)\Gamma (1-s)\Gamma (s)}\over {\Gamma (1+r+n)\Gamma (1-s+r+n)}}&=\left\{\matrix{
\up{4}{{(-1)^n\Gamma (s-n)}\over {n!}}\hfill&\hbox{for }r=0\cr \down{4}{{(-1)^n\Gamma (-s-n)}\over {n!}}\hfill&\hbox{for }r=s\cr} \right\}\cr
\MS&={{(-1)^n}\over {n!}}\,\Gamma (s-2r-n),\cr} \eqno (3.12)$$
where the last equality holds true in view of the fact that $r\in \{0,s \}$.
Combining (3.9) and (3.12), we have found the desired formula for the coefficients in the continuation formula (2.14): $$\eqalign{g_n(r)&={{(-1)^n}\over {n!}}\,\,{{\Gamma (s-2r-n)\,\prod\limits_{ j=1}^{p+1}\,\Gamma (a_j+r+n)}\over {\prod\limits_{j=1}^p\,\Gamma ( b_j+r+n)}}\,\,\lim_{m\to\infty}\,{1\over {\Gamma (s-r-n)}}\cr \MS\cdot\,&_{p+2}F_{p+1}\left(\matrix{\up{4}\hfill a_1+r+n,a_ 2+r+n,\cdots ,a_{p+1}+r+n,-m;\cr
\down{4}b_1+r+n,b_2+r+n,\cdots ,b_p+r+n,1-s+r+n-m;\cr} \,1\right),\cr}
\eqno (3.13)$$
where $r\in \{0,s\}$ and $r$ may be set equal to $s$ only after the limit has been evaluated. The hypergeometric series of unit argument is Saalsch\"utzian and terminating because of the numerator parameter $ -m$.
Such series can be summed for $p=1$ by means of the formula (2.4) of Saalsch\"utz, and for larger $p$ analogous transformation formulas are available ({\it cf}. [7] and [8]). It is by means of these formulas that the limit $m\to\infty$ can be performed and (3.13) becomes meaningful. For example, we may use
$$\eqalign{
&{1\over {\Gamma (b_1)\Gamma (b_2)}}\,\,_4F_3\left(\matrix{ \up{4}\hfill a_1,a_2,a_3,-m;\cr
\down{4}\hfill b_1,b_2,1-s-m;\cr}
\,1\right)\cr
\MS&\qquad ={{(a_1+s)_m(a_2+s)_m(a_3)_m}\over {(s)_m\Gamma (b_1+m) \Gamma (b_2+m)}}\,_4F_3\left(\matrix{
\up{4}\hfill b_1-a_3,b_2-a_3,s,-m;\cr
\down{4}\hfill a_1+s,a_2+s,1-a_3-m;\cr}
\,1\right)\cr
\MS&\quad (s=b_1+b_2-a_1-a_2-a_3\neq 0,-1,-2,\cdots ;\,\,\, 
m=0,1,2,\cdots )\cr}
\eqno (3.14)$$
from [7] in order to obtain
$$\eqalign{
g_n(r)&={{(-1)^n}\over {n!}}\,\,
{{\Gamma (a_1+r+n)\Gamma (a_2+r+n)
\Gamma (s-2r-n)}\over {\Gamma (a_1+s)\Gamma (a_2+s)}}\cr \MS&\qquad\cdot {}_3F_2\left(\matrix{
\up{4}\hfill b_1-a_3,b_2-a_3,s-r-n;\cr
\down{4}\hfill a_1+s,a_2+s;\cr}
\,1\right)}\eqno (3.15)$$
for $p=2$ in accordance with (3.32) in [5] and with (2.10) and (2.15) above. In a similar way, the transformation: $$\eqalign{&{1\over {\Gamma (b_1)\Gamma (b_2)\Gamma (b_3)}}\,\,_5F_ 4\left(\matrix{\up{4}\hfill a_1,a_2,a_3,a_4,-m;\cr \down{4}\hfill b_1,b_2,b_3,1-s-m;\cr}
\,1\right)\cr
\MS&\quad ={{(a_1+s)_m(a_2+s)_m(a_3)_m(a_4)_m}\over {(s)_m\Gamma ( b_1+m)\Gamma (b_2+m)\Gamma (b_3+m)}}\cr
\MS&\qquad\cdot\sum_{k=0}^m\,{{(b_1+b_3-a_3-a_4)_k(b_2+b_3-a_3-a_4 )_k(1-b_3-m)_k(s)_k(-m)_k}\over {(a_1+s)_k(a_2+s)_k(1-a_3-m)_k(1-a_ 4-m)_k\,\,k!}}\cr
\MS&\qquad\cdot\,_4F_3\left(\matrix{\up{4}\hfill b_3-a_3,b_3-a_4,a_ 1+a_2+s+m-1,-k;\cr
\down{4}\hfill b_1+b_3-a_3-a_4,b_2+b_3-a_3-a_4,b_3+m-k;\cr} \,1\right)\cr
\MS&\quad (s=b_1+b_2+b_3-a_1-a_2-a_3-a_4\neq 0,-1,-2,\cdots ;\,\,\, m=0,1,2,\cdots )\cr}
\eqno (3.16)$$
yields
$$\eqalign{g_n(r)&={{(-1)^n}\over {n!}}\,\,{{\Gamma (a_1+r+n)\Gamma (a_2+r+n)\Gamma (s-2r-n)}\over {\Gamma (a_1+s)\,\Gamma (a_2+s)}}\cr \MS&\qquad\cdot\sum_{k=0}^{\infty}\,{{(b_1+b_3-a_3-a_4)_k(b_2+b_3- a_3-a_4)_k(s-r-n)_k}\over {(a_1+s)_k(a_2+s)_k\,\,k!}}\cr \MS&\qquad\cdot\,_3F_2\left(\matrix{\up{4}\hfill b_3-a_3,b_3-a_4,- k;\cr
\down{4}b_1+b_3-a_3-a_4,b_2+b_3-a_3-a_4;\cr} \,1\right)\cr}
\eqno (3.17)$$
for $p=3$. Finally, using the transformation: 
$$\eqalign{
&{1\over {\Gamma (b_1)\Gamma (b_2)\Gamma (b_3)\Gamma (b_4)}}\, 
\,_6F_5\left(\matrix{
\up{4}\hfill a_1,a_2,a_3,a_4,a_5,-m;\cr
\down{4}\hfill b_1,b_2,b_3,b_4,1-s-m;\cr} \,\,1\right)\cr
\MS&\qquad ={{(a_1+s)_m(a_2+s)_m(a_3)_m(a_4)_m(a_5)_m}\over {(s)_m\Gamma (b_1+m)\Gamma (b_2+m)\Gamma (b_3+m)\Gamma (b_4+m)}}\cr \MS&\cdot\sum_{k=0}^m{{(b_1+b_3+b_4-a_3-a_4-a_5)_k(b_2+b_3+b_4-a_ 3-a_4-a_5)_k}\over {(a_1+s)_k(a_2
+s)_k}}\cr
\MS&\qquad\qquad\cdot {{(1-b_3-m)_k
(1-b_4-m)_k(s)_k(-m)_k}\over
{(1-a_3-m)_k(1-a_4-m)_k(1-a_5-m)_k\,\,k!}}\cr \MS&\cdot\sum_{\ell =0}^k{{(b_3+b_4-a_3-a_5)_{\ell} 
(b_3+b_4-a_4-a_5)_{\ell}}\over
{(b_1+b_3+b_4-a_3-a_4-a_5)_{\ell}(b_2+b_3+b_4-a_3-a_4-a_5)_{\ell} }}\cr
\MS&\qquad\qquad\qquad\cdot{{(a_5+m-k)_{\ell}(a_1+a_2+s+m-1)_{\ell}(-k)_{\ell}}\over 
{(b_3+m-k)_{\ell}(b_4+m-k)_{\ell}\,\ell !}}\cr \MS&\cdot\,_4F_3\left(\matrix{
\up{4}\hfill b_3-a_5,b_4-a_5,b_3+b_4-a_3-a_4-a_5-m+k,-\ell ;\cr \down{4}\hfill b_3+b_4-a_3-a_5,b_3+b_4-a_4-a_5,1-a_5-m+k-\ell ;\cr} \,\,1\right)\cr
\MS&(s=b_1+b_2+b_3+b_4-a_1-a_2-a_3-a_4-a_5\neq 0,-1,-2,\cdots ; 
\,\,\,m=0,1,2,\cdots ),\cr}
\eqno (3.18)$$
we can obtain
$$\eqalign{
&g_n(r)={{(-1)^n}\over {n!}}\,{{\Gamma (a_1+r+n)\Gamma 
(a_2+r+n)\Gamma (s-2r-n)}\over {\Gamma (a_1+s)\Gamma (a_2+s)}}\cr \MS&\quad\cdot\sum_{k=0}^{\infty}{{(b_1+b_3+b_4-a_3-a_4-a_5)_k(b_ 
2+b_3+b_4-a_3-a_4-a_5)_k(s-r-n)_k}\over {(a_1+s)_k(a_2+s)_k\,\,k!}}\cr \MS&\quad\cdot\sum_{\ell =0}^k{{(b_3+b_4-a_3-a_5)_{\ell}(b_3+b_4- 
a_4-a_5)_{\ell}(-k)_{\ell}}\over {(b_1+b_3+b_4-a_3-a_4-a_5)_{\ell} (b_2+b_3+b_4-a_3-a_4-a_5)_{\ell}\,\,\ell !}}\cr \MS&\quad\cdot\,_3F_2\left(\matrix{
\up{4}\hfill b_3-a_5,b_4-a_5,-\ell ;\cr
\down{4}\hfill b_3+b_4-a_3-a_5,b_3+b_4-a_4-a_5;\cr} \,\,1\right)\cr}
\eqno (3.19)$$
for $p=4$. Writing these results in the form of (2.15), we may get the corresponding expressions:
$$\eqalign{A_k^{(3)}&={{(b_1+b_3-a_3-a_4)_k(b_2+b_3-a_3-a_4)_k}\over { k!}}\cr
\MS&\qquad\cdot\,_3F_2\left(\matrix{\up{4}\hfill b_3-a_3,b_3-a_4, -k;\cr
\down{4}\hfill b_1+b_3-a_3-a_4,b_2+b_3-a_3-a_4;\cr} \,\,1\right)\cr}
\eqno (3.20)$$
and
$$\eqalign{A_k^{(4)}&={{(b_1+b_3+b_4-a_3-a_4-a_5)_k(b_2+b_3+b_4-a_ 3-a_4-a_5)_k}\over {k!}}\cr
\MS&\quad\cdot\sum_{\ell =0}^k{{(b_3+b_4-a_3-a_5)_{\ell}(b_3+b_4- a_4-a_5)_{\ell}(-k)_{\ell}}\over {(b_1+b_3+b_4-a_3-a_4-a_5)_{\ell} (b_2+b_3+b_4-a_3-a_4-a_5)_{\ell}\,\ell !}}\cr \MS&\quad\cdot\,_3F_2\left(\matrix{\up{4}\hfill b_3-a_5,b_4-a_5,- \ell ;\cr
\down{4}\hfill b_3+b_4-a_3-a_5,b_3+b_4-a_4-a_5;\cr} \,\,1\right).\cr}
\eqno (3.21)$$
These representations are different from (2.11) and (2.12), respectively, and this implies the existence of certain identities considered in [7]. There is
another set of Saalsch\"utzian-like formulas [8], however, such that when they are applied in (3.13) they do yield just (2.11) and (2.12) above.
\bigskip
\bigskip\newline
\noindent{\bf 4.~~Zero-balanced series}
\bigskip
There are two ways to treat the $s$-balanced series, that is, the case when $s$ is equal to an integer. We may start from the formulas for general $s$, Equations (2.14) and (2.15), add $-\epsilon$ to each of the parameters and, as a consequence, replace $s$ by $s+\epsilon$. Assuming then that $s$ is an integer, we have to perform the limit $\epsilon\to 0$. This procedure here is not more complicated than in the Gaussian case, $p=1$, due to the fact that the crucial dependence on $a_1$, $a_2,$ $s$, $r$ remains essentially the same for $p=2,3,\cdots$ or, stated otherwise, the $A_k^{(p)}$ are not modified by the replacements but remain independent of $\epsilon$. Alternatively, if we want to avoid the task of the limit procedure, we may start from (2.9) and use the known analytic continuation formulas of the Gaussian hypergeometric series from Equations (15.3.10) to (15.3.12) of [1] for instance. Both
of these ways lead to the same results displayed in [6] for arbitrary $s$. Here we want to restrict our attention to the simplest case, the zero-balanced series. We have, if $s$ according to (1.2) is equal to zero, $$\eqalign{&{{\Gamma (a_1)\Gamma (a_2)\cdots\Gamma (a_{p+1})}\over { \Gamma (b_1)\cdots\Gamma (b_p)}}\,\,_{p+1}F_p\left(\matrix{\up{4}\hfill a_1,a_2,\cdots ,a_{p+1};\cr
\down{4}\hfill b_1,\cdots ,b_p;\cr}
\,\,z\right)\cr
\MS&\qquad =\sum_{n=0}^{\infty}\,\,{{(a_1)_n(a_2)_n}\over {n!\,\,n !}}\biggl\{\sum_{k=0}^n\,{{(-n)_k}\over {(a_1)_k(a_2)_k}}\,A_k^{(p )}\cr
\MS&\qquad\cdot\left\{\psi (1+n-k)+\psi (1+n)-\psi (a_1+n)-\psi ( a_2+n)-\log (1-z)\right\}\cr
\MS&\qquad +(-1)^n\,(n)!\,\sum_{k=n+1}^{\infty}\,{{(k-n-1)!}\over { (a_1)_k(a_2)_k}}\,A_k^{(p)}\biggr\}(1-z)^n\cr} \eqno (4.1)$$
for $|1-z|<1$, $|\arg (1-z)|<\pi$, and $p=2,3,\cdots$, where the $ A_k^{(p)}$ are given
by (2.10) to (2.13), or (3.20) and (3.21), and the condition $$\Re (a_3+n)>0\land\cdots\land\Re (a_{p+1}+n)>0$$ is required for convergence of the infinite series. 

Restricting our attention to the leading terms, when $z\to 1$, we have $$\eqalign{&{{\Gamma (a_1)\Gamma (a_2)\cdots\Gamma (a_{p+1})}\over { \Gamma (b_1)\cdots\Gamma (b_p)}}\,\,\,_{p+1}F_p\left(\matrix{\up{4}\hfill a_1,a_2,\cdots ,a_{p+1};\cr
\down{4}\hfill b_1,\cdots ,b_p;\cr}
\,\,z\right)\cr
\MS&\qquad =L(p)\left[1+O(1-z)\right]-\log (1-z)\left[1+O(1-z)\right ]\cr}
\eqno (4.2)$$
with
$$L(p)=2\psi (1)-\psi (a_1)-\psi (a_2)+B(p),\eqno (4.3)$$ where, if the condition $\Re (a_j)>0$ for $j=3,4,\cdots ,p+1$ is satisfied, $$B(p)=\sum_{k=1}^{\infty}\,{{(k-1)!}\over {(a_1)_k(a_2)_k}}\,A_k^{ (p)}.\eqno (4.4)$$
While we obtain
$$\eqalign{B(2)&=\sum_{k=1}^{\infty}\,{{(b_2-a_3)_k(b_1-a_3)_k}\over { k(a_1)_k(a_2)_k}}\cr
\MS&={{(b_1-a_3)(b_2-a_3)}\over {a_1a_2}}\,\,_4F_3\left(\matrix{\up{ 4}\hfill b_1-a_3+1,b_2-a_3+1,1,1;\cr
\down{4}\hfill a_1+1,a_2+1,2;\cr}
\,\,1\right),\cr}
\eqno (4.5)$$
in accordance with (1.4), the $B(p)$ for $p\geq 3$ appear in different forms depending on which of the various ways discussed above they are derived. We have
$$\eqalign{B(3)&=\sum_{k=1}^{\infty}\,{{(b_3+b_2-a_4-a_3)_k(b_1-a_ 3)_k}\over {k(a_1)_k(a_2)_k}}\cr
\MS&\qquad\qquad\cdot\,_3F_2\left(\matrix{\up{4}\hfill b_3-a_4,b_ 2-a_4,-k;\cr
\down{4}\hfill b_3+b_2-a_4-a_3,1+a_3-b_1-k;\cr} \,\,1\right),\cr}
\eqno (4.6)$$
if we use (4.4) together with (2.11), or $$\eqalign{B(3)&=\sum_{k=1}^{\infty}\,{{(b_1+b_3-a_3-a_4)_k\,(b_2+ b_3-a_3-a_4)_k}\over {k(a_1)_k(a_2)_k}}\cr \MS&\qquad\qquad\cdot\,_3F_2\left(\matrix{\up{4}\hfill b_3-a_3,b_ 3-a_4,-k;\cr
\down{4}\hfill b_1+b_3-a_3-a_4,b_2+b_3-a_3-a_4;\cr} \,\,1\right)\cr}
\eqno (4.7)$$
if we use (4.4) together with (3.20). Similarly, we find from (4.4) and (2.12) that
$$\eqalign{B(4)&=\sum_{k=1}^{\infty}\,{{(b_4+b_3+b_2-a_5-a_4-a_3)_ k(b_1-a_3)_k}\over {k(a_1)_k(a_2)_k}}\cr \MS&\qquad\cdot\sum_{\ell =0}^k\,{{(b_4+b_3-a_5-a_4)_{\ell}(b_2-a_ 4)_{\ell}(-k)_{\ell}}\over {(b_4+b_3+b_2-a_5-a_4-a_3)_{\ell}(1+a_3 -b_1-k)_{\ell}\,\ell !}}\cr
\MS&\qquad\cdot\,_3F_2\left(\matrix{\up{4}\hfill b_4-a_5,b_3-a_5, -\ell ;\cr
\down{4}\hfill b_4+b_3-a_5-a_4,1+a_4-b_2-\ell ;\cr} \,\,1\right),\cr}
\eqno (4.8)$$
or from (4.4) and (3.21) that
$$\eqalign{&B(4)=\sum_{k=1}^{\infty}{{(b_1+b_3+b_4-a_3-a_4-a_5)_k( b_2+b_3+b_4-a_3-a_4-a_5)_k}\over {k(a_1)_k(a_2)_k}}\cr \MS&\quad\cdot\sum_{\ell =0}^k\,{{(b_3+b_4-a_3-a_5)_{\ell}(b_3+b_ 4-a_4-a_5)_{\ell}(-k)_{\ell}}\over {(b_1+b_3+b_4-a_3-a_4-a_5)_{\ell} (b_2+b_3+b_4-a_3-a_4-a_5)_{\ell}\,\ell !}}\cr \MS&\quad\cdot\,_3F_2\left(\matrix{\up{4}\hfill b_3-a_5,b_4-a_5,- \ell ;\cr
\down{4}\hfill b_3+b_4-a_3-a_5,b_3+b_4-a_4-a_5;\cr} \,\,1\right).\cr}
\eqno (4.9)$$

Alternatively, starting from (2.8) with $p=2$ and using (2.5) and (2.3), we obtain, because of $s=0$, (4.5) above again. In a similar way, with $p=3$ in (2.8) and application of (4.2) and (4.3) with (4.5), we find that
$$\eqalign{
B(3)&={{(b_1-a_3)(b_2+b_3-a_3-a_4)}\over {a_1a_2}}\cr \MS&\qquad\qquad\cdot {}_4F_3\left(\matrix{ 
\up{4}\hfill b_1-a_3+1,b_2+b_3-a_3-a_4+1,1,1;\cr \down{4}\hfill a_1+1,a_2+1,2;\cr}
\,\,1\right)\cr
\MS&\qquad +{{(b_2-a_4)(b_3-a_4)\Gamma (a_1)\Gamma (a_2)\Gamma (a_ 3)}\over {\Gamma (b_1)\Gamma (b_2+b_3-a_4+1)}}\cr \MS&\qquad\cdot\sum_{\ell =0}^{\infty}\,{{(b_2-a_4+1)_{\ell}(b_3- a_4+1)_{\ell}}\over {(\ell +1)!\,(b_2+b_3-a_4+1)_{\ell}}}\,\,_3F_2\left(\matrix{\up{4}\hfill a_1,a_2,a_3;\cr
\down{4}\hfill b_1,b_2+b_3-a_4+1+\ell ;\cr} \,\,1\right).\cr}
\eqno (4.10)$$
The general formula of this type turns out to be $$\eqalign{
B(p)&={{(b_1-a_3)}\over {a_1a_2}}\left(
\sum_{j=2}^p\,b_j-\sum_{j=3}^{p+1}\,a_j\right)\cr \MS&\qquad\qquad\cdot\,_4F_3\left(\matrix{ 
\up{8}\hfill b_1-a_3+1,\sum\limits_{j=2}^p\, 
b_j-\sum\limits_{j=3}^{p+1}\,a_j+1,1,1;\cr \down{8}\hfill a_1+1,\,a_2+1,2;\cr}
\,\,1\right)\cr
\MS&\quad +\sum_{k=3}^p\,{{(b_{k-1}-a_{k+1})\left(\sum\limits_{j=k}^ p\,b_j-\sum\limits_{j=k+1}^{p+1}\,a_j\right)\,\Gamma (a_1)\cdots\Gamma (a_k)}\over {\Gamma (b_1)\cdots\Gamma (b_{k-2})\,\Gamma\left(\sum\limits_{ j=k-1}^p\,b_j-\sum\limits_{j=k+1}^{p+1}\,a_j+1\right)}}\cr \MS&\quad\cdot\sum_{\ell =0}^{\infty}\,{{(b_{k-1}-a_{k+1}+1)_{\ell}\left (\sum\limits_{j=k}^p\,b_j-\sum\limits_{j=k+1}^{p+1}\, a_j+1\right)_{\ell}}\over {(\ell +1)!\,\left(\sum\limits_{j=k-1}^p\,b_j-\sum\limits_{ j=k+1}^{p+1}\,a_j+1\right)_{\ell}}}\cr
\MS&\quad\cdot\,_kF_{k-1}\left(\matrix{\up{6}\hfill a_1,a_2,\cdots ,a_k;\cr
\down{6}\hfill b_1,\cdots ,b_{k-2},\,\sum\limits_{j=k-1}^p\,b_j-\sum\limits_{ j=k+1}^{p+1}\,a_j+1+\ell ;\cr}
\,\,1\right).\cr}
\eqno (4.11)$$
Various other representations can be found. An example for $p=3$ is $$\eqalign{B(3)&=\sum_{k=1}^{\infty}\,{{(b_1-a_3)_k(b_2+b_3-a_3-a_ 4)_k}\over {k(a_1)_k(a_2)_k}}\cr
\MS&\qquad +\sum_{k=1}^{\infty}{{(b_2-a_4)_k(b_3-a_4)_k}\over {k(a_ 1)_k(a_2)_k}}\,\,_3F_2\left(\matrix{\up{4}\hfill b_1-a_3,b_2+b_3-a_ 3-a_4+k,k;\cr
\down{4}\hfill a_1+k,a_2+k;\cr}
\,\,1\right).\cr}
\eqno (4.12)$$
This follows from (2.7) if the term for $m=0$ is evaluated separately by means of (4.2) and (4.3) and if the hypergeometric series in the remaining sum is evaluated at $z=1$ by means of (2.16) and (2.10). It is closely related to (4.10) if the hypergeometric series there in the second term is transformed by means of (2.16). It also follows from (4.6) if there the sum over $k$ and the sum implied by the hypergeometric series are interchanged. Of the various possibilities demonstrated here for $p=3$ as an example, the representation (4.6) or (4.7) has the advantage that only a single infinite series occurs, the hypergeometric series involved consisting of a finite number of terms because of the numerator parameter $-k$.\goodbreak \bigskip
\bigbreak\newline
\noindent{\bf 5.~~Partial sums of zero-balanced series} 
\bigskip
A related problem is to find the asymptotic behaviour of the partial sums of hypergeometric series at unit argument in the zero-balanced case. This problem was solved by Ramanujan ({\it cf}. [15], [3], and [4]) and later proved by Evans and Stanton [9] for $p=2$. More generally, we have, if $s$ according to (1.2) is equal to zero, $$\eqalign{&\sum_{k=0}^{m-1}\,{{\Gamma (a_1+k)\Gamma (a_2+k)\cdots \Gamma (a_{p+1}+k)}\over {\Gamma (b_1+k)\cdots\Gamma (b_p+k)\Gamma (1+k)}}\cr
\MS&\qquad ={\cal L}(p)-\psi (1)+\log (m)+O(m^{-1}),\qquad m\to \infty\cr}
\eqno (5.1)$$
with
$${\cal L}(p)=\sum_{k=0}^{\infty}\left\{{{\Gamma (a_1+k)\Gamma (a_2+k)\cdots\Gamma (a_{p+1}+k)}\over {\Gamma (b_1+k)\cdots\Gamma (b_p+k)\Gamma (1+k)}}-{1\over {k+1}}\right\}.\eqno (5.2)$$ To derive (5.1) we observe that, in the zero-balanced case, the ratio of the gamma functions is equal to $k^{-1}\left\{1+O(k^{-1})\right\}$ and so each summand is $O(k^{-2})$ as $k\to\infty$. Therefore, we have for the partial sum $${\cal L}(p)=\sum_{k=0}^{m-1}\left\{{{\Gamma (a_1+k)\Gamma (a_ 2+k)\cdots\Gamma (a_{p+1}+k)}\over {\Gamma (b_1+k)\cdots\Gamma (b_ p+k)\Gamma (1+k)}}-{1\over {k+1}}\right\}+O(m^{-1})\eqno (5.3)$$ as $m\to\infty$. Also we have
$$\sum_{k=0}^{m-1}\,{1\over {k+1}}=-\psi (1)+\psi (1+m)=-\psi (1)+\log m+\ts m^{-1}+O(m^{-2})\eqno (5.4)$$
as $m\to\infty$, by means of Stirling's formula. Combining (5.3) and (5.4), we get (5.1). On the other hand, from (4.2) we have $$\sum_{k=0}^{\infty}\,{{\Gamma (a_1+k)\Gamma (a_2+k)\cdots\Gamma (a_{p+1}+k)}\over {\Gamma (b_1+k)\cdots\Gamma (b_p+k)\Gamma (1+k)}}\, z^k=L(p)-\log (1-z)+o(1)\eqno (5.5)$$
as $z\to 1$. Expanding the logarithmic term and putting it on the left-hand side, we arrive at
$$\sum_{k=0}^{\infty}\,\left\{{{\Gamma (a_1+k)\Gamma (a_2+k)\cdots \Gamma (a_{p+1}+k)}\over {\Gamma (b_1+k)\cdots\Gamma (b_p+k)\Gamma (1+k)}}-{z\over {k+1}}\right\}\,z^k=L(p)+o(1).\eqno (5.6)$$ The limit as $z\to 1$ now yields
$${\cal L}(p)=L(p),\eqno (5.7)$$
and this is the really interesting point: the constant needed here in (5.1) for the asymptotic behaviour of the partial sum of a zero-balanced hypergeometric series at unit argument is equal to the constant needed above, in (4.2) and (4.3), for the analytic continuation of the infinite series. For example, (5.7) with (4.3) and $p=3$ yields $${\cal L}(3)=2\psi (1)-\psi (a_1)-\psi (a_2)+B(3),\eqno (5.8)$$ where any of the various representations of $B(3)$ according to (4.6), (4.7), (4.10), or (4.12) above may be inserted. The use of (4.10) corresponds to the formula given by A.K. Srivastava [19]. All the other results of the same author [19] concern the partial sums of zero-balanced series with additional relations between the parameters, a topic beyond the scope of the present investigation. \bigskip
\bigskip\newline
\noindent{\bf Acknowledgments}
\bigskip
The present investigation was completed during the second-named author's visit to the University of Heidelberg in July 1997. This work was supported, in part, by the {\it Natural Sciences and Engineering Research Council of Canada } under Grant OGP0007353.
\bigskip
\frenchspacing
\baselineskip=12pt
\parskip 3pt plus 1pt minus 1pt
\noindent {\bf References}
\medskip
\ref M. Abramowitz and I.A. Stegun (Editors), {\it Handbook of Mathematical Functions}, Dover, New York, 1965. 

\ref W.N. Bailey, {\it Generalized Hypergeometric Series}, Stechert-Hafner, New York, 1964.

\ref B.C. Berndt, Chapter 11 of Ramanujan's second notebook, {\it Bull. London Math. Soc.} {\bf 15}(1983), 273--320. 

\ref B.C. Berndt, {\it Ramanujan's Notebooks}, Part II, Springer-Verlag, New York, 1989.

\ref W. B\"uhring, The behavior at unit argument of the hypergeometric function $_3F_2$, {\it SIAM J. Math. Anal.} {\bf 18}(1987), 1227--1234.

\ref W. B\"uhring, Generalized hypergeometric functions at unit argument, {\it Proc. Amer. Math. Soc.} {\bf 114}(1992), 145--153. 

\ref W. B\"uhring, Transformation formulas for terminating Saalsch\"utzian hypergeometric series of unit argument, {\it J. Appl. Math. Stochast}. {\it Anal}. {\bf 8}(1995), 189--194. 

\ref W. B\"uhring, Terminating Saalsch\"utzian hypergeometric series at unit argument, {\it Math. Japon}. {\bf 44}(1996), 15--23. 

\ref R.J. Evans and D. Stanton, Asymptotic formulas for zero-balanced hypergeometric series, {\it SIAM J. Math. Anal.} {\bf 15}(1984), 1010--1020.

\ref R.J. Evans, Ramanujan's second notebook: Asymptotic expansions for hypergeometric series and related functions, in {\it Ramanujan Revisited} (G.E. Andrews {\it et al}., Editors) (Proc. of the Ramanujan Centenary Conference, 1987), Academic Press, New York, 1988, 537--560.

\ref Y.L. Luke, {\it The Special Functions and Their Approximations}, Vol. 1, Academic Press, New York, 1969.

\ref O.I. Marichev and S.L. Kalla, Behaviour of hypergeometric function $_pF_{p-1}(z)$ in the vicinity of unity, {\it Rev. T\'ecn. Fac. Ingr. Univ. Zulia} {\bf 7}(1984), 1--8.

\ref N. N\o rlund, Hypergeometric functions, {\it Acta Math}. {\bf 94}(1955), 289--349.

\ref F.W.J. Olver, {\it Asymptotics and Special Functions}, Academic Press, New York, 1974.

\ref S. Ramanujan, {\it Notebooks}, Vol. 2, Tata Institute of Fundamental Research, Bombay, 1957.

\ref M. Saigo, On properties of the Appell hypergeometric functions $F_2$ and $F_3$ and the generalized Gauss function $_3F_2$, {\it Bull. Central. Res. Inst. Fukuoka Univ.} {\bf 66}(1983), \penalty-10000 27--32. 

\ref M. Saigo and H.M. Srivastava, The behavior of the zero-balanced hypergeometric series $_pF_{p-1}$ near the boundary of its convergence region, {\it Proc. Amer. Math. Soc}. {\bf 110}(1990), 71--76. 

\ref L.J. Slater, {\it Generalized Hypergeometric Functions}, 
Cambridge University Press, \penalty-10000 Cambridge, 1966. 

\ref A.K. Srivastava, Asymptotic behaviour of certain zero-balanced hypergeometric \penalty-10000 series, {\it Proc. Indian Acad. Sci. Math. Sci}. {\bf 106}(1996), 39--51. 

\ref H.M. Srivastava and P.W. Karlsson, {\it Multiple Gaussian Hypergeometric Series}, \penalty-10000 Halsted Press (John Wiley and Sons), New York, 1985.

\bye